\documentclass{article}

\usepackage{amsmath, amsthm, amssymb, amscd}

\newtheorem{definition}{Definition}[section]
\newtheorem{lemma}[definition]{Lemma}
\newtheorem{proposition}[definition]{Proposition}
\newtheorem{example}[definition]{Example}
\newtheorem{theorem}[definition]{Theorem}
\newtheorem{corollary}[definition]{Corollary}
\newtheorem{remark}[definition]{Remark}
\newcommand{\bo}{\mathbf}
\newcommand{\bb}{\mathbb}
\newcommand{\ot}{\otimes}
\renewcommand{\L}{{\mathcal{L}}}
\newcommand{\M}{{\mathcal{M}}}
\newcommand{\bra}{[\hspace{-0.17em}[}
\newcommand{\ket}{]\hspace{-0.17em}]}
\newcommand{\p}{{\prime}}
\newcommand{\ti}{\tilde}

\newcommand{\cc}{\cdot}
\renewcommand{\c}{\circ}
\renewcommand{\l}{{\langle}}
\renewcommand{\r}{{\rangle}}
\renewcommand{\a}{\alpha}
\renewcommand{\b}{\beta}
\newcommand{\e}{\varepsilon}

\newcommand{\A}{\mathcal{A}}

\newcommand{\Hom}{\mathrm{Hom}}
\newcommand{\Der}{\mathrm{Der}}
\newcommand{\End}{\mathrm{End}}

\title{Courant brackets on noncommutative algebras\\
and\\
omni-Lie algebras}

\author{Kyousuke UCHINO}

\date{}

\begin{document}
\maketitle
\begin{abstract}
We define a Courant bracket on an associative
algebra using the theory of Hochschild homology,
and we introduce the notion of Dirac algebra.
We show that the bracket of an omni-Lie algebra is
quite a kind of Courant bracket.
\end{abstract}

\section{Introduction.}
T. Courant \cite{Cou} defines a skew-symmetric bracket (\ref{oricou1}) below
on the set of sections
$\Gamma (TM\oplus T^{*}M)$ on a smooth manifold $M$
\begin{equation}\label{oricou1}
\bra(X,\a),(Y,\b)\ket_{skew}:=
([X,Y],\L_{X}\b-\L_{Y}\a+\frac{1}{2}d(\l Y,\a\r-\l X,\b\r)),
\end{equation}
where $(X,\a),(Y,\b)\in \Gamma(TM\oplus T^{*}M)$.
The bracket is not a Lie bracket, but the modified bracket
\begin{equation}\label{oricou12}
\bra(X,\a),(Y,\b)\ket:=([X,Y],\L_{X}\b-\L_{Y}\a+d\l Y,\a\r)
\end{equation}
satisfies a Leibniz identity and the bracket (\ref{oricou1}) is given as
the skew-symmetrization of (\ref{oricou12}).
These brackets (\ref{oricou1}) and (\ref{oricou12}) are
both called {\em Courant brackets}.
In addition, he gives a smooth nondegenerate symmetric bilinear form
on $TM\oplus T^{*}M$:
\begin{equation}\label{oricou2}
((x,a),(y,b)):=\frac{1}{2}(\l y,a\r+\l x,b\r),
\end{equation}
where $(x,a),(y,b)\in TM\oplus T^{*}M$.
The Courant bracket and the bilinear form are used to give
a characterization of Poisson structure on $M$.
Let $\pi$ be a 2-vector field on $M$, and let $L_{\pi}$
denote the graph of $\pi$, i.e.,
the set of elements $(\ti{\pi}(a),a)$,
where $a\in T^{*}M$ and $\ti{\pi}:T^{*}M\to TM$ is the bundle map
defined by $\pi(a_{1},a_{2})=\l\ti{\pi}(a_{1}),a_{2}\r$.
$\pi$ is a Poisson structure if and only if
the Courant bracket is closed on the set of sections
$\Gamma L_{\pi}$ and $L_{\pi}$ is maximally isotropic
for the bilinear form (\ref{oricou2}).
Such subbundles of $TM\oplus T^{*}M$
are called {\em Dirac structures} (\cite{Cou}).
\begin{definition}\label{coudirac}
Let $M$ be a smooth manifold. A subbundle $L$ of $TM\oplus T^{*}M$
is called a Dirac structure, if the Courant bracket (\ref{oricou12}),
or equivalently (\ref{oricou1})
is closed on the set of sections $\Gamma L$ and $L$ is maximally isotropic
for the bilinear form (\ref{oricou2}).
\end{definition}
A. Weinstein \cite{We} gives a {\em linearization} of (\ref{oricou1}),
or (\ref{oricou12}) motivated by an integrability problem of
Courant brackets.
We refer \cite{Kin} for the study of the integrability problem
of Courant brackets.
Let $V$ be a vector space.
Weinstein's bracket is defined on the space $gl(V)\oplus V$:
\begin{equation}\label{omnibras}
\bra(\xi_{1},v_{1}),(\xi_{2},v_{2})\ket_{skew}:=
([\xi_{1},\xi_{2}],\frac{1}{2}(\xi_{1}(v_{2})-\xi_{2}(v_{1}))),
\end{equation}
where $(\xi_{1},v_{1}),(\xi_{2},v_{2})\in gl(V)\oplus V$.
This bracket is the skew-symmetrization of a Leibniz bracket:
\begin{equation}\label{omnibral}
\bra(\xi_{1},v_{1}),(\xi_{2},v_{2})\ket:=([\xi_{1},\xi_{2}],\xi_{1}(v_{2})).
\end{equation}
The $V$-valued nondegenerate symmetric
bilinear form is also defined by
\begin{equation}\label{omnibi}
((\xi_{1},v_{1}),(\xi_{2},v_{2}))=\frac{1}{2}(\xi_{2}(v_{1})+\xi_{1}(v_{2})).
\end{equation}
Similar to Poisson structures on a manifold, every Lie algebra structure
on $V$ is characterized as the graph.
Let $\mu:V\ot V\to V$ be a binary operation.
Set the graph of $\mu$: $L_{\mu}:=\{(\ti{\mu}(v),v) \ | \ v\in V\}$,
where $\ti{\mu}:V\to gl(V)$ is the map defined by $\ti{\mu}(v)(u)=\mu(v,u)$.
The operation $\mu$ is a Lie bracket if and only if Weinstein's bracket
(\ref{omnibras}), or equivalently (\ref{omnibral})
is closed on $L_{\mu}$ and $L_{\mu}$ is maximally isotropic
for the bilinear form (\ref{omnibi}). Such objects are called
{\em D-structures} in \cite{We}.
He calls $gl(V)\oplus V$ an {\em omni-Lie algebra}.
Here we consider relationships between Courant brackets (\ref{oricou12})
and Weinstein's brackets (\ref{omnibras}), or (\ref{omnibral}).\\
\indent
In \cite{We} it is suggested that $V$ is a non-unital algebra
of linear functions on the dual space $V^{*}$ with trivial multiplication.
Then $gl(V)$ is the set of derivations of $V$.
Furthermore $(0,v)\in gl(V)\oplus V$ is
a certain derivative $D:v\mapsto (0,v)$, similar to
manifolds cases $D:C^{\infty}(M)\to\Gamma (TM\oplus T^{*}M)$,
$f\mapsto (0,df)$. So one can view omni-Lie algebras
as geometrical ``linearization'' of Courant's original examples.
\medskip\\
\indent
In this paper we construct an algebraic Courant bracket
using Hochschild cohomology (resp. homology) groups.
Let $\A$ be an associative and unital algebra,
not necessarily commutative, and we set the Hochschild
cohomology (resp. homology) group $H^{1}(\A,\A)$ (resp. $H_{1}(\A,\A)$).
In Section 3, we define on the space $H^{1}(\A,\A)\oplus H_{1}(\A,\A)$
a Leibniz bracket by the same formula as (\ref{oricou12}),
using algebraic derivatives.
We will call the bracket on $H^{1}(\A,\A)\oplus H_{1}(\A,\A)$
a \textbf{Courant bracket} on $\A$.
Denote $H^{1}(\A,\A)\oplus H_{1}(\A,\A)$ by $E(\A)$.
Our motivation is given by the following example.
\begin{example}\label{motiex}
Let $V$ be a finite-dimensional vector space.
Set $V[1]:=V\oplus \bb{R}\cc 1$ as a unital algebra over the field $\bb{R}$,
where the multiplication is almost trivial except the unit $1$.
Then $H^{1}(V[1],V[1])$ is just $gl(V)$ and the Courant bracket on $V[1]$
has the same formula as Weinstein's bracket (\ref{omnibral}):
\[
\bra(\xi_{1},dv_{1}),(\xi_{2},dv_{2})\ket=
([\xi_{1},\xi_{2}],d\xi_{1}(v_{2})),
\]
where $d:V[1]\to H_{1}(V[1],V[1])$ is an algebraic de Rham derivative.
(See Section 5, for the detailed study.)
\end{example}
In addition, the symmetric bilinear form of $E(\A)$ is also well-defined
by the formula (\ref{oricou2}) without the factor $1/2$,
by means of a duality between $H^{1}(\A,\A)$ and $H_{1}(\A,\A)$.
We wish a {\em nondegenerate} symmetric bilinear form to define
the notion of Dirac structure on $\A$. However, the bilinear form
on $E(\A)$ is degenerate in general.
We notice that the kernel of the bilinear form becomes an ideal
for the Courant bracket on $\A$. Thus we have the exact sequence
of Leibniz algebras:
\begin{equation}
\begin{CD}
0@>>>J@>>>E(\A)@>>>E(\A)/J@>>>0,
\end{CD}
\end{equation}
where $J$ is the set of the kernel of the bilinear form.
The quotient Leibniz algebra $\e(\A):=E(\A)/J$ has an induced nondegenerate
symmetric bilinear form and the induced Courant bracket.
Even if $\A$ is noncommutative, thanks to the nondegeneracy of
the bilinear form on $\e(\A)$, the notion of \textbf{Dirac structure}
is well-defined as a maximally isotropic submodule $L$ of $\e(\A)$
such that the induced Courant bracket is closed on the submodule.
We call the pair $(\A,L)$ a (noncommutative) \textbf{Dirac algebra}.
We will show that every Poisson bracket
on a commutative algebra is characterized as the corresponding
Dirac structure.\\
\indent
We denote the matrix algebra of an algebra $\A$ by $\M_{r}(\A)$.
In Proposition \ref{morita}
we will show that the Courant bracket on $\M_{r}(\A)$
is isomorphic to the one of $E(\A)$
and the bilinear form is also preserved by the isomorphism.
By this proposition, we obtain an Courant bracket
isomorphism $\e(\A)\cong \e(\M_{r}(\A))$.
The first main theorem of this paper is
\begin{theorem}\label{morita2005}
Let $\A$ be a unital and associative algebra.
Then there exists a Courant bracket
isomorphism $\e(\A)\cong \e(\M_{r}(\A))$
preserving the bilinear form.
Thus Dirac structures on $\A$ and $\M_{r}(\A)$ correspond bijectivelly.
\end{theorem}
\indent
It is well-known that the dual bundle of a Lie algebroid
$A\to M$ is a Poisson manifold with Lie-Poisson bracket.
When $M$ is a point, the Lie algebroid is a Lie algebra
and the Lie-Poisson bracket is the ordinary one.
One can view the algebra $V[1]$ of Example \ref{motiex}
as a linearization of the smooth functions on the vector bundle
$V^{*}\to \{o\}$ on a point. In fact the part $\bb{R}\cc 1$ is the set of
functions on the base point. Thus $\e(V[1])$ is the linearization of
Courant's original type example $TV^{*}\oplus T^{*}V^{*}$.
The second main result of this paper is
\begin{theorem}\label{main}
Let $V$ be a vector space of finite dimension.
Then $\e(V[1])$ is isomorphic to omni-Lie algebra $gl(V)\oplus V$,
i.e., Weinstein's bracket on $gl(V)\oplus V$
is the (induced) Courant bracket on $\e(V[1])$.
\end{theorem}

The paper organized as follows.\\
\indent
In Section 2 we recall some basic properties of Hochschild (co)homology
theory and the algebraic operations corresponding to Lie derivative,
interior product and exterior derivative.\\
\indent
In Section 3 we define the Courant bracket, the bilinear form on
$E(\A)$ and study the basic property.
Especially we show that the algebraic Courant bracket on $\A$
satisfies the axioms of Courant algebroids.\\
\indent
In subsection 4.1 we study the bilinear form and introduce
the quotient space $\e(\A)$ with nondegenerate bilinear form.
Algebraic Dirac structures are introduced (Definition \ref{coudirac}).
Theorem \ref{morita2005} is proved.\\
\indent
In subsection 4.2 we show that a Poisson algebra is a Dirac algebra
and every Poisson bracket is characterized by the corresponding
Dirac structure.\\
\indent
In Section 5 the second main theorem is shown.\\

\noindent {\bf Acknowledgements.}
I would like to thank very much Professors
Yoshiaki Maeda and Akira Yoshioka for their helpful comments
and encouragement.
\section{Preliminalies}
In this section we recall Hochschild (co)homology groups of algebras
and set an algebraic differential-calculus.
We refer the book \cite{Lo} for the detailed study of the theory
of Hochschild (co)homology.
\subsubsection{Hochschild homology}
Let $k$ be a commutative ring, $\A$ be an algebra over
the ring $k$. The Hochschild $n$-complex is
$C_{n}(\A,\A):=\A\otimes \A^{\otimes n}$,
where the tensor product is defined over $k$.
The boundary map $b:C_{n}(\A,\A)\to C_{n-1}(\A,\A)$ is defined by
the rule below.
Let $P_{i}:C_{n}(\A,\A)\to C_{n-1}(\A,\A)$ be a $k$-homomorphism:
\begin{eqnarray*}
P_{i}(a_{0}\ot...\ot a_{n})&:=&
(-1)^{i}(a_{0}\ot...\ot a_{i}a_{i+1}\ot...\ot a_{n}), (0\le i\le n-1)\\
P_{n}(a_{0}\ot...\ot a_{n})&:=&(-1)^{n}(a_{n}a_{0}\ot...\ot a_{n-1}),
\end{eqnarray*}
where $a_{0},...,a_{n}\in \A$. The map $b$ is defined by the formula:
\begin{equation*}
b(a_{0}\ot...\ot a_{n}):=\sum^{n}_{i=0}P_{i}(a_{0}\ot...\ot a_{n}).
\end{equation*}
It holds that $b^{2}=0$, and thus the homology groups
$H_{n}(\A,\A)$ are defined.
For example, since $b(a_{0}\ot a_{1})=[a_{0},a_{1}]=a_{0}a_{1}-a_{1}a_{0}$,
the $0$-th Hochschild homology group is $H_{0}(\A,\A)=\A/[\A,\A]$,
where $[\A,\A]$ is a $k$-module generated by
all $[a,a^{\p}]$. We denote the center of $\A$ by $Z(\A)$.
One can check that by the action of $Z(\A)\to C_{n}(\A,\A)$:
$z(a_{0}\ot...\ot a_{n})=(za_{0}\ot...\ot a_{n})$,
each $H_{n}(\A,\A)$ becomes a $Z(\A)$-module.
In fact, for any $z\in Z(\A)$ we obtain
$zb(a_{0}\ot...\ot a_{n})=b(za_{0}\ot...\ot a_{n})$.
If $\A$ is commutative then $H_{0}(\A,\A)=\A$, and if $\A$ is unital then
$H_{1}(\A,\A)$ is isomorphic to the $\A$-module of
{\em K\"{a}hler differentials} which is an $\A$-module
generated by 1-forms $ada^{\p}$ (see the next subsection 2.0.2.).
\subsubsection{K\"{a}hler differentials}
We assume $\A$ is unital and commutative. Set an $\A$-module $O_{\A|k}$
generated by $da$ for any $a\in\A$, where $d$ is merely a symbol.
Define two relations (or axioms) on the module $O_{\A|k}$:
\begin{eqnarray}
d(\lambda a+\lambda^{\p}a^{\p})-\lambda da
-\lambda^{\p}da^{\p}&=&0,\label{defk1}\\
d(aa^{\p})-ada^{\p}-a^{\p}da&=&0,\label{defk2}
\end{eqnarray}
where $\lambda,\lambda^{\p}\in k$.
The quotient module $O_{\A|k}/\sim$ is the module
of K\"{a}hler differentials, and it is denoted by $\Omega^{1}_{\A|k}$.
It is known that $H_{1}(\A,\A)\cong \Omega^{1}_{\A|k}$
(see 1.1.10 Proposition in \cite{Lo}).
The isomorphism between $H_{1}(\A,\A)$ and $\Omega^{1}_{\A|k}$ is given
by $a_{0}\ot a_{1}\cong a_{0}da_{1}$, on the level of cycles.
In fact, by the relation (\ref{defk1}), $O_{\A|k}$ becomes
the tensor product $\A\ot \A$, and the second relation is
the same as the defining relation of the Hochschild homology
$H_{1}(\A,\A)$.
\subsubsection{Hochschild cohomology}
Next, we consider the Hochschild cohomology groups for general algebras.
The $n$-complex $C^{n}(\A,\A)$
is $\Hom_{k}(\A^{\ot n},\A)$ and when $n=0$, $C^{0}(\A,\A)=\A$.
The coboundary map $\b$ is defined by the following formula.
For any $f\in C^{n}(\A,\A)$:
\begin{multline*}
\b(f)(a_{1}\ot...\ot a_{n+1})=a_{1}f(a_{2}\ot...\ot a_{n+1})\\
+\sum^{n}_{i=1} (-1)^{i}f(a_{1}\ot...\ot a_{i}a_{i+1}\ot...\ot a_{n+1})\\
+(-1)^{n+1}f(a_{1}\ot...\ot a_{n})a_{n+1}
\end{multline*}
and $\b:C^{0}(\A,\A)\to C^{1}(\A,\A)$ is $a\mapsto [a,\cdot]$
for any $a\in\A$. It is easily checked that $H^{0}(\A,\A)=Z(\A)$
and the cocycles of $C^{1}(\A,\A)$ is the set of derivations on $\A$.
Denote the derivations on $\A$ by $\Der(\A)$. We have
\[
H^{1}(\A,\A)=\Der(\A)/[\A,\cc],
\]
where $[\A,\cc]$ is the submodule of $C^{1}(\A,\A)$ generated by
inner derivations $[a,\cc]:a^{\p}\mapsto [a,a^{\p}]$.
Especially if $\A$ is commutative then $H^{1}(\A,\A)=\Der(\A)$.
Note that each $H^{n}(\A,\A)$ is also a $Z(\A)$-module.
\subsubsection{Algebraic derivatives}
Secondly, we recall a Lie bracket on $H^{1}(\A,\A)$, a Lie derivative
$\L_{X}$, an interior product $i_{X}$ and Connes' boundary map $B$
on homology.
\begin{remark}
In \cite{Lo}, the Lie derivative and the interior product are
denoted by $L_{D}$ and $e_{D}$ respectively.
Here we use geometrical notations $\L_{X}$ and $i_{X}$.
\end{remark}
\indent
For any $\A$, $\Der(\A)$ has a canonical Lie bracket by
taking the commutator. One can easily check that the module generated
by inner-derivations $[a,\cdot]$ is an ideal of $\Der(\A)$.
Thus a Lie bracket on $H^{1}(\A,\A)$ is induced.\\
\indent
A Lie derivative $\L_{X}:H_{n}(\A,\A)\to H_{n}(\A,\A)$ associated
with an element $X\in \Der(\A)$ is defined by the formula below, on the level
of cycles (Section 4.1 of \cite{Lo}).
\begin{equation*}
\L_{X}(a_{0}\ot...\ot a_{n}):=\sum^{n}_{i=0}a_{0}
\ot...\ot X(a_{i})\ot...\ot a_{n}.
\end{equation*}
Further, we can show that $\L_{[a,\cdot]}$ is the zero map on the level of
homology groups for any $a\in\A$ (see 4.1.5 Proposition in \cite{Lo}).
Thus $\L_{X}$ is well-defined for any $X\in H^{1}(\A,\A)$.\\
\indent
It is known that an interior product is also well-defined
on Hochschild homology. For any $X\in \Der(\A)$, set an operator
$i_{X}:H_{n}(\A,\A)\to H_{n-1}(\A,\A)$,
on the level of cycles, by the formula:
\begin{equation*}
i_{X}(a_{0}\ot...\ot a_{n}):=(-1)^{n+1}(X(a_{n})a_{0}\ot...\ot a_{n-1}).
\end{equation*}
\begin{lemma}
For any $a\in\A$,
$i_{[a,\cc]}$ is the zero map on the level of homology groups.
\end{lemma}
\begin{proof}
For any $a^{\p}\in\A$,
set the map $h_{a^{\p}}:C_{n}(\A,\A)\to C_{n}(\A,\A)$,
$a_{0}\ot...\ot a_{n}\mapsto a^{\p}a_{0}\ot...\ot a_{n}$.
We have
\begin{eqnarray*}
h_{a^{\p}}\c b(a_{0}\ot...\ot a_{n})&=&a^{\p}a_{n}a_{0}\ot...\ot a_{n-1}+
\sum^{n-1}_{i=0} P_{i}(a^{\p}a_{0}\ot...\ot a_{n}),\\
b\c h_{a^{\p}}(a_{0}\ot...\ot a_{n})&=&a_{n}a^{\p}a_{0}\ot...\ot a_{n-1}+
\sum^{n-1}_{i=0} P_{i}(a^{\p}a_{0}\ot...\ot a_{n}).
\end{eqnarray*}
Thus $h_{a^{\p}}\c b-b\c h_{a^{\p}}=(-1)^{n+1}i_{[a^{\p},\cc]}$
which implies that $i_{[a^{\p},\cc]}$ is homotopic to the zero map.
\end{proof}
Thus the interior product $i_{X}$ is well-defined for any
$X\in H^{1}(\A,\A)$. Like smooth manifold cases,
the following lemma holds.
\begin{lemma}\label{lem2}
For any $X,Y\in H^{1}(\A,\A)$:
\[
\L_{[X,Y]}=\L_{X}\c\L_{Y}-\L_{Y}\c\L_{X},\ \
i_{[X,Y]}=\L_{X}\c i_{Y}-i_{Y}\c \L_{X}.
\]
\end{lemma}
\begin{proof}
We only show the second formula.
For the first formula, we refer 4.1.6 Corollary in \cite{Lo}.
For any $\bo{a}:=a_{0}\ot...\ot a_{n}$:
\begin{multline}\label{lem21}
(-1)^{n+1}\L_{X}\c i_{Y}(\bo{a})=\L_{X}Y(a_{n})a_{0}\ot...\ot a_{n-1}\\
=XY(a_{n})a_{0}\ot...\ot a_{n-1}+Y(a_{n})X(a_{0})\ot...\ot a_{n-1}+\\
\sum^{n-1}_{i=1} Y(a_{n})a_{0}\ot...\ot X(a_{i})\ot...\ot a_{n-1},
\end{multline}
and on the other hand,
\begin{multline}\label{lem22}
(-1)^{n+1}i_{Y}\c \L_{X}(\bo{a})=(-1)^{n+1}i_{Y}\sum^{n}_{i=0} a_{0}\ot...\ot X(a_{i})\ot...\ot a_{n}=\\
Y(a_{n})X(a_{0})\ot...\ot a_{n-1}+
\sum^{n-1}_{i=1}Y(a_{n})a_{0}\ot...\ot X(a_{i})\ot...\ot a_{n-1}+\\
YX(a_{n})a_{0}\ot...\ot a_{n-1}.
\end{multline}
The difference of (\ref{lem21}),(\ref{lem22}) is $(-1)^{n+1}i_{[X,Y]}(\bo{a})$.
\end{proof}
\subsubsection{Connes' boundary map}
In the following, we assume that $\A$ is unital.
The Connes' boundary map $B:H_{n}(\A,\A)\to H_{n+1}(\A,\A)$ is defined by
using cyclic operators on $C_{n}(\A,\A)$ (Section 2.1.7 of \cite{Lo}).
We use an explicit definition:
\begin{multline*}
B(a_{0}\ot...\ot a_{n}):=\sum^{n}_{i=0}(-1)^{ni}(1\ot a_{i}\ot...\ot a_{n}\ot a_{0}\ot...\ot a_{i-1})+\\
(-1)^{ni}(a_{i}\ot1\ot a_{i+1}\ot...\ot a_{n}\ot a_{0}\ot...\ot a_{i-1}),
\end{multline*}
where $1$ is the unit of $\A$.
For example,
\[
B(a_{0}\ot a_{1})=1\ot a_{0}\ot a_{1}-1\ot a_{1}\ot a_{0}+
a_{0}\ot 1\ot a_{1}-a_{1}\ot 1\ot a_{0}.
\]
\begin{remark}
It is known that the condition of boundary operator $B^{2}=0$
is satisfied. However, in our explicit definition,
it is difficult to show the condition.
\end{remark}

It is known that $\L_{X}=B\c i_{X}+i_{X}\c B$ for each $H_{n}(\A,\A)$
(4.1.9 Corollary of \cite{Lo} and see Remark \ref{norm} below).
We directly show the condition:
$\L_{X}=B\c i_{X}+i_{X}\c B$ for $H_{1}(\A,\A)$.
For any cycles $\a$ and $\a^{\p}$
we denote $\a \equiv \a^{\p}$, if $\a=\a^{\p}$
on the level of homology.
\begin{lemma}\label{lem1}
For any $X\in \Der(\A)$, and any cycle $\a\in C_{1}(\A,\A)$:
\[
\L_{X}(\a)\equiv B\c i_{X}(\a)+i_{X}\c B(\a).
\]
Thus $\L_{X}=B\c i_{X}+i_{X}\c B$ on $H_{1}(\A,\A)$
for any $X\in H^{1}(\A,\A)$.
\end{lemma}
\begin{proof}
We can put $\a=a_{0}\ot a_{1}$ without loss of generality.
By $b(1\ot 1\ot a)=a\ot 1$, $a\ot 1\equiv 0$.
Thus we obtain
\begin{multline*}
i_{X}\c B(a_{0}\ot a_{1})=i_{X}(1\ot a_{0}\ot a_{1}-1\ot a_{1}\ot a_{0}+
a_{0}\ot 1\ot a_{1}-a_{1}\ot 1\ot a_{0})\\
\equiv -X(a_{1})\ot a_{0}+X(a_{0})\ot a_{1},
\end{multline*}
and $B\c i_{X}(a_{0}\otimes a_{1})\equiv 1\ot X(a_{1})a_{0}$.
In addition, we have
\[
b(1\ot X(a_{1})\ot a_{0})=X_{1}(a_{1})\ot a_{0}-1\ot X(a_{1})a_{0}
+a_{0}\ot X(a_{1}).
\]
Thus $1\ot X(a_{1})a_{0}\equiv X_{1}(a_{1})\ot a_{0}+a_{0}\ot X(a_{1})$.
This gives a proof of the lemma.
\end{proof}
Here we remark that when $\A$ is commutative, the derivative $da$ into
the space of K\"{a}hler differentials is the same as
the the boundary $B(a)$:
\[
B(a)\equiv 1\ot a \cong da,\ \ \ B=d:\A\to \Omega^{1}_{\A|k}.
\]
\begin{remark}\label{norm}
We can take a normalized-Hochschild homology group
$\overline{H}_{n}(\A,\A)$ which is defined by the certain
quotient $C_{n}(\A,\A)/\sim$ of Hochschild complex.
It is known that the normalized-Hochschild
homology group is isomorphic with an ordinary one
(see 1.1.14 of \cite{Lo}).
The condition $\L_{X}=B\c i_{X}+i_{X}\c B$ is shown
in the normalized framework for any $H_{n}(\A,\A)$.
\end{remark}
\subsubsection{Canonical pairings}
We now set the pairing between $H^{1}(\A,\A)$ and $H_{1}(\A,\A)$
using the interior product by the form:
\begin{equation}\label{pairing}
\l\cc,\cc\r:H^{1}(\A,\A)\times H_{1}(\A,\A)\to H_{0}(\A,\A),\ \
\l X,\a\r:=i_{X}\a,
\end{equation}
where $X\in H^{1}(\A,\A)$ and $\a\in H_{1}(\A,\A)$.
Note that the pairing is $Z(\A)$-bilinear.\\
\indent
We remark here that the pairing (\ref{pairing}) is equivalent to
the {\em Kronecker product}. The Kronecker product
$\l\cc,\cc\r:H^{n}(\A,\A)\ot H_{n}(\A,\A)\to \A\ot_{\A^{e}}\A$,
is a canonical pairing between cohomology groups and homology
groups defined by, on the level of (co)chains,
\[
\l f,a_{0}\ot a_{1}\ot...\ot a_{n}\r=f(a_{1}\ot...\ot a_{n})\ot_{A^{e}}a_{0},
\]
where $f\in C^{n}(\A,\A)$, $\A^{e}:=\A\ot \A^{op}$ and $\A^{op}$
is the opposite algebra of $\A$ (see 1.5.9 Duality of \cite{Lo}).
One can easily show that $\A\ot_{\A^{e}}\A \cong H_{0}(\A,\A)$.
The isomorphism is $a\ot_{A^{e}}a^{\p}\cong\overline{aa^{\p}}$,
where $\overline{aa^{\p}}$ is the equivalence class of $aa^{\p}$.
In fact, by the definition,
we have $a\ot_{A^{e}}a^{\p}=1(1\ot a)\ot_{A^{e}}a^{\p}=1\ot a^{\p}a$.
On the other hand,
$a\ot_{A^{e}}a^{\p}=1(a\ot 1)\ot_{A^{e}}a^{\p}=1\ot aa^{\p}$.
Thus $a\ot_{A^{e}}a^{\p}=a^{\p}\ot_{A^{e}}a$.
This commutativity is expressed as the {\em abelianzation}
$\A/[\A,\A]=H_{0}(\A,\A)$.
\medskip\\
\indent
Recall the bilinear forms (\ref{oricou2}) and (\ref{omnibi}).
By means of the bilinear form, the notion of Dirac structure
is defined as a maximally isotropic subspace.
We use carefully the term ``maximally isotropic"
in the algebraic framework.
\medskip\\
\indent
Let $k$ be a unital commutative ring, and let $E$ and
$M$ be (left) $k$-modules, and let $(\cdot,\cdot)$ be a $M$-valued
nondegenerate symmetric $k$-bilinear form on $E$.
Here $(\cdot,\cdot)$ is nondegenerate, namely $(e,\cc): E\to M$
is injective for any nontrivial $e\in E$.
\begin{definition}\label{defmaximal}
Under the notations above,
let $L$ be a submodule of $E$. We say that $L$ is ``isotropic"
for the bilinear form, if the bilinear form is zero on $L$.
When $L$ is isotropic, we say that $L$ is ``maximally isotropic",
if $(e,\cc)$ vanishes on $L$ then $e$ is in $L$ for any $e\in E$.
\end{definition}
\section{Courant bracket of $H^{1}(\A,\A)\oplus H_{1}(\A,\A)$}
In this section we define a Courant bracket on an associative algebra
using the operations of Section 2.
\begin{definition}
Let $\A$ be a unital and associative $k$-algebra.
We call a bracket on $H^{1}(\A,\A)\oplus H_{1}(\A,\A)$ below
a Courant bracket on $\A$.
\[
\bra (X_{1},\a_{1}),(X_{2},\a_{2})\ket=
([X_{1},X_{2}],\L_{X_{1}}\a_{2}-\L_{X_{2}}\a_{1}+
B\l X_{2},\a_{1}\r),
\]
where $(X_{1},\a_{1}),(X_{2},\a_{2})\in H^{1}(\A,\A)\oplus H_{1}(\A,\A)$.
We denote $H^{1}(\A,\A)\oplus H_{1}(\A,\A)$ by $E(\A)$.
\end{definition}
When $k$ contains $1/2$,
we have Courant's original formula (\ref{oricou1})
as the skew-symmetrization of the Courant bracket on $\A$.
\medskip\\
\indent
We set a symmetric $Z(A)$-bilinear form $(\cdot,\cdot)$ on
$H^{1}(\A,\A)\oplus H_{1}(\A,\A)$ using the formula
(\ref{oricou2}) without the factor $1/2$, i.e.,
for any $e_{1}:=(X_{1},\a_{1}),e_{2}:=(X_{2},\a_{2})\in E(\A)$:
\begin{equation}\label{bilinearform1219}
(e_{1},e_{2}):=\l X_{2},\a_{1}\r+\l X_{1},\a_{2}\r.
\end{equation}
Note that this bilinear form is $H_{0}(\A,\A)$-valued in general.
In addition, we set a map $\rho:E(\A)\to H^{1}(\A,\A)$
as the canonical projection:
\begin{equation}
\rho(X,\a):=X.
\end{equation}
By the definition, $\rho$ has a $Z(A)$-linearity.
We notice a derivative action:
\begin{equation}\label{act}
H^{1}(\A,\A)\times Z(\A)\to Z(\A),\ \ (X,z)\mapsto X(z),
\end{equation}
This action is well-defined on the level of homology,
since $[a,z]=0$ for any $a\in\A$.
\begin{proposition}\label{prop1}
Let $\A$ be a unital and associative $k$-algebra.
Then the Courant bracket satisfies the following properties.
For any $e_{1},e_{2},e_{3}\in E(\A)$ and $z\in Z(\A)$:
\begin{eqnarray}
\bra e_{1},\bra e_{2},e_{3}\ket\ket&=&\bra\bra e_{1},e_{2}\ket,e_{3}\ket+
\bra e_{2},\bra e_{1},e_{3}\ket\ket\label{c0}\\
\rho\bra e_{1},e_{2}\ket&=&[\rho(e_{1}),\rho(e_{2})],\label{c1}\\
\bra e_{1},ze_{2}\ket&=&z\bra e_{1},e_{2}\ket+\rho(e_{1})(z)e_{2},\label{c2}\\
2\bra e_{1},e_{1}\ket&=&D(e_{1},e_{1}),\label{c4}\\
\L_{\rho(e_{1})}(e_{2},e_{3})&=&(\bra e_{1},e_{2}\ket,e_{3})+(e_{2}, \bra e_{1},e_{3}\ket),\label{c3}
\end{eqnarray}
where
$D$ is a $k$-homomorphism:
\[
D:H_{0}(\A,\A)\to H_{1}(\A,\A),\ \ \a\mapsto (0,B(\a)),
\]
and $\rho(e_{1})(z)$ of (\ref{c2}) is the action (\ref{act}).
\end{proposition}
\begin{proof}
The formulas (\ref{c1}) and (\ref{c4}) are clearly.
For (\ref{c0}), (\ref{c2}) and (\ref{c3}), by Lemma \ref{lem2}, \ref{lem1}
we can take the same proof as the case of $TM\oplus T^{*}M$
on a smooth manifold.
\end{proof}
The conditions (\ref{c0})-(\ref{c3}) above are the set of axioms
of {\em Courant algebroids} in \cite{Roy}
(see also \cite{Kos}, \cite{Liu}).
However $E(\A)$ is not a Courant algebroid, because the bilinear
form is degenerate in general. In the next section we will study
the bilinear form on $E(\A)$.
\medskip\\
\indent
For given algebras $\A$ and $\A^{\p}$, we write $E(\A)\cong E(\A^{\p})$,
if there exists an isomorphism $\phi:H_{0}(\A,\A)\cong H_{0}(\A^{\p},\A^{\p})$
and if there exists a Courant bracket isomorphism preserving the bilinear form
up to $\phi$. We study isomorphisms between Courant brackets.\\
\indent
It is well-known that a unital algebra $\A$ and the matrix algebra
$\M_{r}(\A)$ are Morita equivalent, and thus the Hochschild
(co)homology groups of $\A$ and $\M_{r}(\A)$
are isomorphic (see 1.2.4 and 1.5.6 in \cite{Lo}).
\begin{proposition}\label{morita}
For any $\A$, $E(\A)\cong E(\M_{r}(\A))$.
\end{proposition}
\begin{proof}
We take isomorphisms $cotr:H^{1}(\A,\A)\to H^{1}(\M_{r}(\A),\M_{r}(\A))$
and $inc: H_{1}(\A,\A)\to H_{1}(\M_{r}(\A),\M_{r}(\A))$ in \cite{Lo}.
Here these maps are defined by
\[
cotr(X)(m_{ij}):=(X(m_{ij})),\ \ \
inc(a_{0}\ot...\ot a_{n})=E_{11}(a_{0})\ot...\ot E_{11}(a_{n}),
\]
on the level of chains, where $X\in H^{1}(\A,\A)$, $m_{ij}\in\M_{r}(\A)$
and $E_{11}(a)$ is a matrix such that the $(1,1)$-position is $a$
and other positions are all zero. We denote $cotr$ and $inc$ by $T$
and $I$ respectively.\\
\indent
It is obvious that $T$ is a Lie algebra isomorphism.
First we show that $T\oplus I$ preserves the bilinear form.
It is sufficient to show that
$i_{T(X)}\c I(a_{0}\ot a_{1})=I\c i_{X}(a_{0}\ot a_{1})$.
\begin{eqnarray*}
i_{T(X)}\c I(a_{0}\ot a_{1})&=&i_{T(X)}(E_{11}(a_{0})\ot E_{11}(a_{1}))\\
&=&T(X)(E_{11}(a_{1}))\cc E_{11}(a_{0})\\
&=&E_{11}(X(a_{1}))\cc E_{11}(a_{0})\\
&=&E_{11}(X(a_{1})a_{0})=I\c i_{X}(a_{0}\ot a_{1}).
\end{eqnarray*}
Thus the bilinear form is preserved by the isomorphism.
Secondly we show that the Courant bracket is preserved.
For any $a\in\A$, we have
$B\c I(a)\equiv 1_{\M_{r}(\A)}\ot E_{11}(a)$ and
$I\c B(a)\equiv E_{11}(1)\ot E_{11}(a)$,
where $1_{\M_{r}(\A)}$ is the unit element of $\M_{r}(\A)$.
On the other hand,
\begin{eqnarray*}
(B\c I-I\c B)(a)&\equiv&(1_{\M_{r}(\A)}-E_{11}(1))\ot E_{11}(a)\\
&=&-b\{(1_{\M_{r}(\A)}-E_{11}(1))\ot E_{11}(a)\ot E_{11}(1)\}\\
&\equiv&0.
\end{eqnarray*}
Thus $B\c I(a)=I\c B(a)$ on the level of homology.
We now obtain below, on the level of homology:
$I\c B\c i_{X}(a_{0}\ot a_{1})=B\c i_{T(X)}\c I(a_{0}\ot a_{1})$.
One can directly show:
$i_{T(X)}\c B\c I(a_{0}\ot a_{1})=I\c i_{X}\c B(a_{0}\ot a_{1})$,
on the level of homology.
Thus we obtain $\L_{T(X)}\c I(\a)=I\c \L_{X}(\a)$
and $I\c B\l X,\a\r=B\l T(X),I(\a)\r$ for any $X\in H^{1}(\A,\A)$
and $\a\in H_{1}(\A,\A)$. Thus $T\oplus I$ preserves the Courant bracket.
\end{proof}
This proposition will be used to give a proof
of Theorem \ref{morita2005} in the next section.
As an example of other isomorphisms
we can easily check that $E(\A)\cong E(\A^{op})$,
where $\A^{op}$ is the opposite algebra of $\A$.
(We refer E.2.1.4 of \cite{Lo}.)
\begin{example}\label{opiso}
$E(\A)\cong E(\A^{op})$.
\end{example}


\section{Dirac algebras and Poisson brackets}
\subsection{Dirac structures}
Let $M$ be a smooth manifold.
Dirac structures $L$ on $M$ are defined
as maximally isotropic subbundles of $TM\oplus T^{*}M$
for the bilinear form (\ref{oricou2})
such that the Courant bracket (\ref{oricou12})
is closed on the set of sections $\Gamma L$.
The maximality condition is well-defined
because the bilinear form (\ref{oricou2}) is nondegenerate.
In Courant's original example, the pair
$(M,L)$ is called a Dirac manifold. In this subsection, we
introduce a notion of {\em Dirac algebra}. First we study
the bilinear form of $E(\A)$.
\medskip\\
\indent
Let $\A$ be a unital $k$-algebra. The bilinear form
(\ref{bilinearform1219}) of $E(\A)$
is degenerate in general. But we can show that the kernel
of the bilinear form is an ideal of $E(\A)$ with respect to
the Courant bracket. Denote the kernel by $J$, i.e.,
\begin{equation}\label{JJ}
J:=\{e\in E(\A) \ | \ (e,e^{\p})=0 \ \text{for any} \ e^{\p}\in E(\A)\}.
\end{equation}

\begin{lemma}\label{kernel}
The kernel $J$ is an ideal of $E(\A)$.
\end{lemma}
\begin{proof}
For any $e\in J$, $e_{1},e_{2}\in E(\A)$,
by (\ref{c3}) in Proposition \ref{prop1}
we have
\[
\L_{\rho(e_{1})}(e,e_{2})=(\bra e_{1},e\ket,e_{2})+(e, \bra e_{1},e_{2}\ket).
\]
Since $e$ is in the kernel, we have $(\bra e_{1},e\ket,e_{2})=0$.
By the definition of Courant bracket,
the skew-symmetrization is
\begin{equation}\label{mD}
\bra e_{1},e\ket-\bra e,e_{1}\ket=2\bra e_{1},e\ket-D(e_{1},e),
\end{equation}
where $D$ was defined in Proposition \ref{prop1}.
This implies that $J$ is a two-side ideal.
\end{proof}
From this lemma, when $J\neq E(\A)$, we obtain a nontrivial
Leibniz algebra $E(\A)/J$ with $H_{0}(\A,\A)$-valued nondegenerate
symmetric bilinear form $(\cc,\cc)$. Here the bilinear form on $E(\A)/J$
is $Z(\A)$-bilinear. We denote $E(\A)/J$ by $\e(\A)$.
When $\e(\A)\neq 0$,
we obtain a Leibniz algebra $\e(\A)$ with a nondegenerate bilinear form
and a (induced) Courant bracket. So we define Dirac structures
on noncommutaitve algebras.
\begin{definition}
Let $A$ be a unital and associative $k$-algebra.
We assume that $\e(\A)\neq 0$.
We call a submodule $L$ of $\e(\A)$ a \textbf{Dirac structure} on $\A$,
if $L$ is maximally isotropic for the induced bilinear form on $\e(\A)$
and the induced Courant bracket on $\e(\A)$ is closed on $L$.
We call the pair $(\A,L)$ a (noncommutative) \textbf{Dirac algebra}.
\end{definition}
In general we have no hope of defining
the map $\rho:\e(\A)\to H^{1}(\A,\A)$.
When $\A$ is commutative, the kernel
$J$ becomes a submodule of $H_{1}(\A,\A)$,
thus the map $\rho$ is well-defined (see Lemma \ref{lemmaJ} below).
\medskip\\
\indent
An algebraic meaning of Dirac structure is that
it is a Lie algebra. By the isotropy of Dirac structure
we have a corollary below.
\begin{corollary}
A Dirac structure $L$ on $k$-algebra $\A$ is a $k$-Lie algebra
and the inverse image $p^{-1}(L)$ of the canonical projection
$p:E(\A)\to \e(\A)$ satisfies the defining conditions of Lie algebroids.
For any $l_{1},l_{2}\in L$ and $z\in Z(\A)$: 
\[
\sigma\bra l_{1},l_{2}\ket=[\sigma(l_{1}),\sigma(l_{2})],\ \ \
\bra l_{1},zl_{2}\ket=z\bra l_{1},l_{2}\ket+\sigma(l_{1})(z)l_{2},
\]
where $\sigma$ is an anchor map defined by the composition
$p^{-1}(L)\overset{\rho}{\to} H^{1}(\A,\A)\to \Der(Z(\A))$
and $[\cc, \cc]$ is a commutator on $\Der(Z(\A))$.
\end{corollary}
\begin{proof}
It is obvious that the Courant bracket on $\A$ is closed
on $p^{-1}(L)$. The anchor map is well-defined by the action
(\ref{act}). Two conditions above follow from (\ref{c1}),
(\ref{c2}) in Proposition \ref{prop1}.
\end{proof}
Note that the above $\sigma$ differs from $\rho$
in Proposition \ref{prop1}.
Especially when $k=\bb{R}$ and $Z(\A)$ is the algebra of
smooth functions on a manifold $M$, $p^{-1}(L)$ is just the space
of sections of a Lie algebroid on $M$.\\
\indent
In the next subsection we will show that a Poisson algebra is
a Dirac algebra with the corresponding Dirac structure.
It is well-known that closed 2-forms on a manifold
define Dirac structures (see \cite{Cou}). Similar to manifold cases,
we obtain a proposition below.
\begin{proposition}\label{lomega}
We assume $\e(\A)\neq 0$.
Let $\omega\in H_{2}(\A,\A)$ be a closed 2-form
in the sense of $B(\omega)=0$ and $i_{X}i_{Y}\omega=-i_{Y}i_{X}\omega$.
Then $p(L_{\omega})$ is a Dirac structure, where $L_{\omega}$
is the set of elements $(X,i_{X}\omega)$ and $p$ is the
canonical projection $p:E(\A)\to \e(\A)$.
\end{proposition}
\begin{proof}
It is obvious that $L_{\omega}$ is isotropic on $E(\A)$ and $\e(\A)$.
By the same way as geometrical cases in \cite{Cou},
one can easily check that the Courant bracket is closed
on $L_{\omega}$. We show that $p(L_{\omega})$ is maximally isotropic.
Recall Definition \ref{defmaximal}.
For any $(X,i_{X}\omega)\in L_{\omega}$,
we assume $((X,i_{X}\omega),(Y,\a))=0$ on $E(\A)$.
Then we have $i_{X}i_{Y}\omega=i_{X}\a$ for any $X$,
thus $(0,i_{Y}\omega-\a)$ is in the kernel $J$.
Thus in $\e(\A)$ we have $\a=i_{Y}\omega$, i.e.,
$p(L_{\omega})$ is maximally isotropic.
\end{proof}
From the proposition above, when $\omega$ is trivial,
the projection of $H^{1}(\A,\A)$ is a Dirac structure.
\medskip\\
\indent
We now give a proof of {\bf Theorem \ref{morita2005}} in Introduction.
\begin{proof}
Using the isomorphism $T\oplus I:E(\A)\cong E(\M_{r}(\A))$
in Proposition \ref{morita}, we obtain a Courant bracket isomorphism
\[
p\c(T\oplus I)\c p^{-1}:\e(\A)\cong \e(\M_{r}(\A))
\]
which preserves the bilinear form on $\e(\A)$ up to the
isomorphism $H_{0}(\A,\A)\cong H_{0}(\M_{r}(\A),\M_{r}(\A))$.
Thus Dirac structures correspond bijectively between $\A$ and $M_{r}(\A)$.
This gives the proof of Theorem \ref{morita2005}.
\end{proof}
From Example \ref{opiso} and Theorem \ref{morita2005},
we obtain $\e(\A)\cong \e(\M_{r}(\A)^{op})$.
Using the theorem we give an example of Dirac algebra
on a smooth manifold.
\begin{example}
Set $\A:=C^{\infty}(M)$ which is the set of smooth functions
on a smooth manifold $M$. Then $\M_{r}(\A)$ is identified
with $\Gamma\End(\bb{R}^{r}\times M)$ which is the space of
sections of the endmorphism bundle of the trivial bundle.
Using the identification $\e(C^{\infty}(M))\cong \Gamma(TM\oplus T^{*}M)$,
we obtain Dirac structures on the algebra $\Gamma\End(\bb{R}^{r}\times M)$
from geometrical (i.e. ordinary) Dirac structures on the manifold.
For instance, for a Poisson structure $\pi$ on $M$,
the graph $L_{\pi}$ is a Dirac structure on $C^{\infty}(M)$.
We can denote the derivation $T(X)\in\Der(\M_{r}(\A))$ for
$X\in\Gamma TM$ in the matrix form
$\left(
\begin{array}{cc}
X & 0 \\
0 & X
\end{array}
\right)$,
where we put $r=2$.
On the other hand, $I(fdg)$ is
$\left(
\begin{array}{cc}
f & 0 \\
0 & 0
\end{array}
\right)
\ot
\left(
\begin{array}{cc}
g & 0 \\
0 & 0
\end{array}
\right)$.
Thus the Dirac structure $p\circ(T\oplus I)\circ p^{-1}(\Gamma L_{\pi})$
has the form, on the level of chains,
\[
\left\{
\left.
\left(
\left(
\begin{array}{cc}
fX_{g} & 0 \\
0 & fX_{g}
\end{array}
\right),
\left(
\begin{array}{cc}
f & 0 \\
0 & 0
\end{array}
\right)
\ot
\left(
\begin{array}{cc}
g & 0 \\
0 & 0
\end{array}
\right)
\right)
\right| \ X\in \Gamma TM, \ f,g\in C^{\infty}(M)
\right\},
\]
where $X_{f}$ is the Hamilton vector field of $f$.
\end{example}

\begin{remark}
Lemma \ref{kernel} is important for Courant algebroids in Poisson geometry.
Given a ``week''-Courant algebroid with degenerate symmetric bilinear form,
we can take the quotient bundle with the induced Leibniz bracket and
the nondegenerate bilinear form. Conversely, we expect that
every Courant algebroid is given in this way.
\end{remark}
\subsection{Poisson algebras}
The purpose of this subsection is to show
that every Poisson bracket is characterized as the corresponding
Dirac structure.
\medskip\\
\textbf{Standing Assumptions.}
We assume $\A$ is a unital and commutative $k$-algebra.
Thus $H^{1}(\A,\A)=\Der(\A)$, $H_{1}(\A,\A)=\Omega^{1}_{\A|k}$,
$H^{0}(\A,\A)=Z(\A)=\A=H_{0}(\A,\A)$ and the boundary map
$B=d:\A\to \Omega^{1}_{\A|k}$.
In addition we assume that $\e(\A)\neq 0$.
The condition is satisfied if there exists a nontrivial derivation on $\A$.
Thus this assumption is always satsified in Poisson geometry.
\medskip\\
\indent
It is known that the derivative
$d:\A\to \Omega^{1}_{\A|k}$ has the universality below
(see 1.3.7-1.3.9 of \cite{Lo}).
For any derivative $\delta:\A\to M$ to an $\A$-module,
there exists a unique map $\phi:\Omega^{1}_{\A|k}\to M$ such that
$\delta=\phi\circ d$, here $\phi$ is $\A$-linear.
\begin{remark}
Usually, the universal derivation of an algebra
is defined as the derivation $d:\A\to I/I^{2}$, where
$I$ is a (non-symmetric) $\A$-bimodule generated by
$1\ot a-a\ot 1$ for any $a\in \A$ and $I/I^{2}$
is the symmetrization of $I$.
One can check that $I/I^{2}\cong \Omega^{1}_{\A|k}$.
\end{remark}

\begin{lemma}\label{lemmaJ}
If $(X,\a)$ is an element of the kernel $J$ (\ref{JJ})
of the bilinear form then $X=0$.
\end{lemma}
\begin{proof}
By the assumption, for any $a\in \A$ we have $((0,da),(X,\a))=0$.
When $da\neq 0$, this gives $X(a)=0$.
Even if $da=0$, by the universality above, we have $X(a)=0$.
\end{proof}
By this lemma, when $\A$ is commutative,
the map $\rho:\e(\A)\to \Der(\A)$ is
induced from $\rho$ on $E(\A)$.
In this case, all conditions (\ref{c0})-(\ref{c3})
of Proposition \ref{prop1} are satisfied on $\e(\A)$.
Thus for a commutative algebra $\A$,
$\e(\A)$ can be viewed as an example of {\em Courant algebra}.
In fact if $k$ includes $1/2$ and $\A$
is commutative then $\e(\A)$ becomes an example of $(k,\A)$ C-algebra.
In \cite{We} an algebraic edition of Courant algebroids is
defined on a non-unital commutative algebra, this is called
a {\em C-algebra}. It was shown that omni-Lie algebra $gl(V)\oplus V$
is a C-algebra on the algebra $V$ with trivial
multiplication. In the next section, we will show that the brackets
of omni-Lie algebras are given by the purely algebraic Courant brackets.
\medskip\\
\indent
Now we consider Poisson algebras (on commutative algebras).
In Poisson manifold cases,
it is well-known that a Poisson bracket $\{\cc,\cc\}$ on $C^{\infty}(M)$ is
equivalent with the Poisson structure $\pi\in\Gamma\bigwedge^{2}TM$
using the definition $\{f,g\}=\pi(df,dg)$ for any $f,g\in C^{\infty}(M)$.
Recall that the Poisson condition $[\pi,\pi]=0$ is equivalent to
the Jacobi law of the bracket $\{\cc,\cc\}$.
The Poisson structure $\pi$ is identified with the bundle map
$\ti{\pi}:T^{*}M\to TM$ by the canonical pairing
$\pi(df,dg)=\l\ti{\pi}(df),dg\r$,  and thus the Poisson bracket
is identified with the Dirac structure $L_{\pi}$ given by the graph
of $\ti{\pi}$. For an arbitrary Poisson algebra $\A$ these
identifications are not always defined. But we can get the Dirac structure
of a Poisson algebra.\\
\indent
Let $\{\cc,\cc\}$ be a $k$-bilinear biderivation on a $k$-algebra $\A$,
not necessarily Poisson bracket. From the universality above,
a Hamiltonianzation $\A\to \Der(\A)$, $a\mapsto \{a,\cc\}$ is
given by the formula: $\{a,\cc\}=\ti{\pi}(da)$ using the unique map
$\ti{\pi}:\Omega^{1}_{\A|k}\to \Der(\A)$. So we obtain the
graph of the map $\ti{\pi}$, which we denote by $L_{\pi}$:
\[
L_{\pi}:=\{(\ti{\pi}(\a),\a) \ | \ \a\in \Omega^{1}_{\A|k}\}.
\]
Note that $L_{\pi}$ is a $\A$-submodule of $E(\A)$.
\begin{proposition}\label{pidirac}
Let $\{\cc,\cc\}$ be a $k$-bilinear biderivation on $\A$,
and we put the corresponding map $\ti{\pi}$.
The bracket is a Poisson bracket if and only if
the pairing $(\cc,\cc)$ on $E(\A)$ is zero on $L_{\pi}$ and
the Courant bracket is closed on $L_{\pi}$.
\end{proposition}
\begin{proof}
We assume that $L_{\pi}$ is isotropic and the Courant bracket
is closed on $L_{\pi}$. For any elements
$(\ti{\pi}(da_{1}), da_{1}),(\ti{\pi}(da_{2}), da_{2})\in L_{\pi}$,
by the isotropy condition,
we have $\l\ti{\pi}(da_{2}),da_{1}\r=-\l\ti{\pi}(da_{1}),da_{2}\r$.
Here
$\l\ti{\pi}(da_{2}),da_{1}\r=i_{\ti{\pi}(da_{2})}(da_{1})=\{a_{2},a_{1}\}$.
This gives the skewsymmetry of the bracket. The Courant bracket
of $(\ti{\pi}(da_{1}), da_{1})$ and $(\ti{\pi}(da_{2}), da_{2})$
has the form:
$([\ti{\pi}(da_{1}),\ti{\pi}(da_{2})],\overline{1\ot\{a_{1},a_{2}\}})$,
here $\overline{1\ot\{a_{1},a_{2}\}}$
is the equivalence class of $1\ot\{a_{1},a_{2}\}$, i.e.,
$\overline{1\ot\{a_{1},a_{2}\}}=d\{a_{1},a_{2}\}$ on $\Omega^{1}_{\A|k}$.
Then we have
\begin{equation}\label{jacobi}
\ti{\pi}(d\{a_{1},a_{2}\})=\{\{a_{1},a_{2}\},\cc\}=
[\ti{\pi}(da_{1}),\ti{\pi}(da_{2})],
\end{equation}
this implies that $\{\cc,\cc\}$ is a Poisson bracket.\\
\indent
Conversely, we assume that $\{\cc,\cc\}$ is a Poisson bracket.
Then we have (\ref{jacobi}) by the
Jacobi identity, i.e., generators of $L_{\pi}$ is closed under the
Courant bracket. Since $\Omega^{1}_{\A|k}$ is generated by
$\{da |a\in \A\}$ as $\A$-module, by (\ref{c2}) in Proposition \ref{prop1}
the Courant bracket is closed on $L_{\pi}$.
The isotropy condition of $L_{\pi}$ is equivalent
to the skewsymmetry of $\{\cc,\cc\}$.
\end{proof}

\begin{lemma}\label{isotropy}
The submodule $L_{\pi}$ of Proposition \ref{pidirac}
is maximally isotropic on $E(\A)$, hence $p(L_{\pi})$
is maximally isotropic on $\e(\A)$,
where $p:E(\A)\to \e(\A)$ is the canonical projection.
\end{lemma}
\begin{proof}
For some element $(X,b^{\p}db)$ in $E(\A)$,
we assume that $((X,b^{\p}db),\cc)=0$ on $L_{\pi}$.
Then for any $a\in \A$,
$((X,b^{\p}db),(\ti{\pi}da,da))=X(a)+b^{\p}\{a,b\}=0$.
When $da\neq 0$, this implies that $X(a)=\ti{\pi}(b^{\p}db)(a)$.
Even if $da=0$, by the universality we obtain
$X(a)=\ti{\pi}(b^{\p}db)(a)=0$.
Thus $X=\ti{\pi}(b^{\p}db)$ which gives that $L_{\pi}$
is maximally isotropic on $E(\A)$.
This implies that $p(L_{\pi})$
is maximally isotropic in $\e(\A)$.
\end{proof}

Here we obtain the main result of this subsection.
\begin{proposition}\label{pLpi}
Let $\{\cc,\cc\}$ be a binary and biderivation on $\A$.
The bracket is a Poisson bracket if and only if
$p(L_{\pi})$ is a Dirac structure,
where $L_{\pi}$ is the same as $L_{\pi}$
in Proposition \ref{pidirac}.
\end{proposition}
\begin{proof}
By $p^{-1}(p(L_{\pi}))=L_{\pi}$.
\end{proof}
Note that since $\A$ is commutative, $\e(\A)$ is identified with
$\Der(\A)\oplus (\Omega^{1}_{\A|k}/J)$. Thus $p(L_{\pi})$ is still
the graph of the induced map.\\
\indent
By Proposition \ref{pLpi}, every Poisson
bracket on a unital commutative algebra $\A$ is characterized
by the Dirac structure of $\e(\A)$.
\begin{remark}
There exists the case  $\e(\A)=0$, for example $\A=k$.
For this case we may always take the trivial Poisson bracket on $\A$.
But this zero Poisson bracket is not characterized
by Dirac structures. This is the difficulty of the algebraic
formulation.
\end{remark}
In Section 5 an example of $\e(\A)$
will be given and studied.
\medskip\\
\noindent
\textbf{Noncommutaive Poisson algebras}.
Finally at this subsection, we consider Poisson structures
associated with Poisson brackets.
Let $\{\cc,\cc\}$ be a Poisson bracket on $\A$. Then we have
$k$-homomorphism $\pi:\A\ot\A\to \A$ by $\pi(a\ot a^{\p})=\{a,a^{\p}\}$.
Since $\{\cc,\cc\}$ is a biderivation,
one can easily check that $\pi$ is a Hochscild
2-cocycle, thus there exists the equivalence class $\pi\in H^{2}(\A,\A)$.
We do not know whether the class satisfies the Poisson condition $[\pi,\pi]=0$
on $H^{3}(\A,\A)$ under the Gerstenhaber bracket.
P. Xu \cite{Xu} showed the converse in noncommutative algebra cases.
If $\Pi\in H^{2}(\A,\A)$ satisfies the Poisson
condition then the center $Z(\A)$ becomes a Poisson
algebra by the bracket $\{z,z^{\p}\}:=[z,[\Pi,z^{\p}]]$.
We do not know whether a noncommutative Poisson structure
$\Pi$ defines the Dirac structure or not, in general.
Here we consider a particular case.
If the matrix algebra $\M_{r}(\A)$ of a commutative algebra $\A$
has a Poisson structure $\Pi$ then $Z(\M_{r}(\A))\cong \A$
is a Poisson algebra, and thus we have a Dirac structure
$L_{\pi}$ on $\A$. By Theorem \ref{morita2005}
we obtain the corresponding Dirac structure on $\M_{r}(\A)$.

\section{Omni-Lie algebras v.s. $\e(\A)$}
In this subsection we will show Theorem \ref{main} in Introduction.
\medskip\\
\indent
Let $V$ be a vector space over the field $\bb{R}$.
Set the vector bundle $V^{*}\to \{o\}$ over a point,
where $V^{*}$ is the dual space of $V$.
The fiber-linearized functions on the bundle
is a vector space $V[1]:=V\oplus \bb{R}\cc 1$
with almost trivial multiplication:
\[
v_{1}\cc v_{2}=0 \ \text{and} \ v\cc 1=1\cc v=v \ \text{for any} \
v_{1},v_{2},v\in V.
\]
It is clear that $V[1]$ is commutative.
\medskip\\
\indent
We now move on to the proof of Theorem \ref{main}.
To show the theorem, we determine the module of
K\"{a}hler differentials $\Omega^{1}_{V[1]|\bb{R}}$.
Recall the definition (\ref{defk1}) and (\ref{defk2}) in Section 2.
Let $\{v_{i}\}$ be a basis of $V$. First we consider the module
$O_{V[1]|\bb{R}}$. It is generated by all elements
$\{d1,vd1,dv,vdv^{\p} \ | \ v,v^{\p}\in V\}$
as infinite $\bb{R}$-module.
Since the multiplication on $V$ is trivial,
$v^{\p\p}vdv^{\p}=0$. First, we assume the linearity (\ref{defk1}) of $d$.
Then the dimension of $O_{V[1]|\bb{R}}$ is reduced to
in $1+2\dim V+(\dim V)^{2}$ and the induced module is generated by
$\{d1,v_{i}d1,dv_{i},v_{i}dv_{j}\}$,
because every element of $V[1]$ is generated by $\{1,v_{i}\}$.
Remark that the dimension is the same as the one of tensor product
$V[1]\ot V[1]$. Secondly, we assume the derivation property
(\ref{defk2}) of $d$. All defining relations are generated by $d1=0$,
$v_{i}d1=0$ and $v_{i}dv_{j}=-v_{j}dv_{i}$.
Here we used $d(v_{i}v_{j})=0$.
Thus the dimension of $\Omega_{V[1]|\bb{R}}$ is $\dim V+_{\dim V}C_{2}$
and it is generated by $\{dv_{i}, {v_{i}dv_{j}}_{(i<j)}\}$ as $\bb{R}$-module,
where $_{\cc}C_{\cc}$ is the combination. Thus we have $\bb{R}$-isomorphism
$\Omega^{1}_{V[1]|\bb{R}}\cong V\oplus \bigwedge^{2}V$ by
$dv\cong v$, $vdv^{\p}=v\wedge v^{\p}$.\\
\indent
It is easy to determine the space of derivations $\Der(V[1])$.
For any $X\in\Der(V[1])$, if $X(v)=v^{\p}+r\cc 1$ for any $v,v^{\p}\in V$
then $0=X(v^{2})=2rv$. Thus we have $r=0$, i.e., $\Der(V[1])\subset gl(V)$.
On the other hand $gl(V)$ becomes the space of derivations of $V[1]$
by the rule $\xi(1):=0$ for any $\xi\in gl(V)$. Thus we obtain
$\Der(V[1])\cong gl(V)$ and
\begin{proposition}\label{EV1}
$E(V[1])\cong gl(V)\oplus V\oplus \bigwedge^{2}V$.
\end{proposition}

Now we compute the Courant bracket on $V[1]$.
For any $(\xi_{1},dv_{1}),(\xi_{2},dv_{2})$, the bracket
has Weinstein's formula (\ref{omnibral}):
\begin{eqnarray*}
\bra(\xi_{1},dv_{1}),(\xi_{2},dv_{2})\ket&=&
([\xi_{1},\xi_{2}],\L_{\xi_{1}}dv_{2}-\L_{\xi_{2}}dv_{1}+d\l\xi_{2},dv_{1}\r)\\
&=&([\xi_{1},\xi_{2}],d(\xi_{1}(v_{2}))),
\end{eqnarray*}
where $(\xi_{1},dv_{1}),(\xi_{2},vdv_{2})\in \e(V[1])$
and the Connes boundary map $B$ is $d$.
For $(\xi_{1},dv_{1}),(\xi_{2},vdv_{2})$, the bracket is
$([\xi_{1},\xi_{2}],0)$ from the triviality of the multiplication.
\begin{lemma}\label{finall}
The kernel $J$ of the bilinear form of $E(V[1])$ is
generated by $\{(0,v_{i}dv_{j})\}$, i.e., $J\cong \bigwedge^{2}V$
\end{lemma}
\begin{proof}
By $i_{\xi}(v_{i}dv_{j})=\xi(v_{j})v_{i}=0$ and Lemma \ref{lemmaJ}.
\end{proof}
From the above lemma, we obtain an isomorphism between
Leibniz algebras:
\[
\e(V[1])\cong gl(V)\oplus V,\ \ \ (\xi,dv)\cong (\xi,v).
\]
One can easily check that by the isomorphism the bilinear
forms are isomorphic. Thus Theorem \ref{main} is proved.
We easily obtain the corollary of the theorem.
\begin{corollary}
A Poisson bracket on $V[1]$ corresponds bijectively
with the Lie bracket on $V$. Thus a Lie-Poisson bracket
on $V^{*}$ corresponds bijectively to the Poisson bracket on $V[1]$.
\end{corollary}
\begin{proof}
By the isomorphism, the Dirac structure of a Poisson bracket
on $V[1]$ corresponds to the graph of a Lie algebra structure on $V$.
\end{proof}
%
%
%
%
%
%

\noindent Kyousuke Uchino, \\
Department of Mathematics,\\
Science University of Tokyo,\\
Wakamiya 26, Shinjyku-ku, Tokyo, 162-0827, Japan\\
e-mail; j1103701@ed.kagu.tus.ac.jp

\end{document}